\documentclass[11pt]{amsart}
\usepackage[margin=1.1in]{geometry}
\usepackage[dvipsnames]{xcolor}
\usepackage{amsmath,amssymb,amsthm,mathtools,enumitem,booktabs,hyperref}
\hypersetup{colorlinks=true,linkcolor=blue!50!black,citecolor=blue!50!black,urlcolor=blue!50!black}

\newtheorem{theorem}{Theorem}[section]
\newtheorem{lemma}[theorem]{Lemma}
\newtheorem{proposition}[theorem]{Proposition}
\newtheorem{corollary}[theorem]{Corollary}

\theoremstyle{definition}
\newtheorem{definition}[theorem]{Definition}
\newtheorem{example}[theorem]{Example}
\theoremstyle{remark}
\newtheorem{remark}[theorem]{Remark}

\newcommand{\E}{\mathbb{E}}
\newcommand{\cP}{\mathcal{P}}
\newcommand{\cQ}{\mathcal{Q}}
\newcommand{\cI}{\mathcal{I}}
\newcommand{\cut}[1]{\|#1\|_\square}
\DeclareMathOperator{\Var}{Var}

\title[An $L\log L$ theory of sparse graph convergence]{An $L\log L$ theory of sparse graph convergence}
\author{Behzad Mehrdad}
\email{behzad.mehrdad86@gmail.com}

\begin{document}

\begin{abstract}
Borgs, Chayes, Cohn and Zhao developed a limit theory for sparse graphs based on $L^p$ upper regularity with $p>1$, in which H\"older's inequality controls the contribution of small sets. We develop the analogous theory in the Orlicz class $L\log L$, with H\"older replaced by the Young--Fenchel inequality for the pair $\Psi(u)=(1+u)\log(1+u)-u$, $\Psi^*(v)=e^v-1-v$. We prove a weak regularity lemma for upper $\Psi$-regular graphons with $4^{O(1/\varepsilon^2)}$ parts, obtained from an energy increment in which the modular $\int\Psi$ itself serves as the energy; that every upper $\Psi$-regular sequence has a subsequence converging in the cut metric to a graphon $W$ with $\int\Psi(W)$ bounded by the regularity constant; that convergent sequences with $L\log L$ limits are upper $\Psi$-regular; and that $L\log L$ balls are compact in the cut metric. The class strictly contains every $L^p$ class with $p>1$, and we give sparse $W$-random graph sequences that are upper $\Psi$-regular but not upper $L^p$ regular for any $p>1$. For graphs with a given degree sequence we prove, by Stein's method of exchangeable pairs with switchings, that when the average degree tends to infinity, the density tends to zero and the degree profile converges in $L^1$ to $\varphi$, uniformly random graphs converge in the normalized cut metric to the rank-one graphon $\varphi\otimes\varphi$, with no condition on the maximum degree; the sequence is upper $\Psi$-regular when $\varphi\in L\log L$. We also show that the maximum-entropy $\beta$-model is entrywise the Chung--Lu model when $d_{\max}^2=o(M)$, and that cut-closeness to the $\beta$-model follows from a lower bound $e^{-o(M)}$ on the probability that it has the prescribed degrees; by a companion counting result such a bound holds at every density with $M\gg n\log n\log(M/n)$. Finally, we show that $\Psi$ is the large deviation rate function of sparse Erd\H{o}s--R\'enyi graphs up to a factor $1+O(1/\log\frac1p)$.
\end{abstract}

\maketitle

\section{Introduction}

\subsection{Background}
Dense graph sequences have limits, graphons, and the theory rests on the weak regularity lemma of Frieze and Kannan \cite{FK99,LS06,BCLSV08,Lov12}. For sparse graphs, with edge density $\rho_n\to0$, one normalizes by $\rho_n$ and must control ``dense spots''. Borgs, Chayes, Cohn and Zhao \cite{BCCZ19,BCCZ18} introduced \emph{$L^p$ upper regularity}: after normalization, every stepping of the graph over a partition with parts of measure at least $\eta$ has $L^p$ norm at most $C$. For $p>1$ they proved a regularity lemma, existence of limits in $L^p$, compactness of $L^p$ balls in the cut metric, and convergence of sparse $W$-random graphs; the case $p=1$ is treated under an additional uniformity hypothesis \cite[Appendix~C]{BCCZ19}. In \cite{BCCZ19} the exponent $p>1$ enters through H\"older's inequality: a set $E$ carries at most $\|U\|_p\,\mu(E)^{1-1/p}$ of the mass of $U$, and $\int_{U>K}U\le\|U\|_p^pK^{1-p}$.

\subsection{This paper}
We replace H\"older by the Young--Fenchel inequality
\begin{equation}\label{eq:young}
uv\;\le\;\Psi(u)+\Psi^*(v)\qquad(u,v\ge0),\qquad \Psi(u)=(1+u)\log(1+u)-u,\quad\Psi^*(v)=e^v-1-v,
\end{equation}
and the $L^p$ norm by the modular $\Lambda(U)=\int\Psi(|U|)$. Since $\Psi(u)\sim u\log u$, the resulting class is $L\log L$, which strictly contains every $L^p$ with $p>1$. The H\"older estimates become (Lemmas~\ref{lem:smallset} and~\ref{lem:tail})
\[
\int_EU\;\le\;\frac{1+\Lambda(U)}{\log(1/\mu(E))},
\qquad
\int_{\{U>K\}}U\;\le\;\frac{\Lambda(U)}{\log(1+K)-1}\quad(K>e-1).
\]

\subsection{Results}
Graphons here are nonnegative (Section~\ref{sec:prelim}).
\begin{itemize}[leftmargin=2em]
\item Theorem~\ref{thm:regularity}: a weak regularity lemma for $(C,\eta)$-upper $\Psi$-regular graphons, with at most $4^N$ parts, $N=\lfloor8C(1+2\Psi^{-1}(C))/\varepsilon^2\rfloor+1$. The energy is the modular $\Lambda(W_\cP)$ itself, so no truncation is needed.
\item Theorem~\ref{thm:limits}: every $C$-upper $\Psi$-regular sequence has a subsequence converging in $\delta_\square$ to a graphon $W$ with $\Lambda(W)\le C$.
\item Proposition~\ref{prop:converse}: if $\delta_\square(W_n,W)\to0$ and $\Lambda(W)<\infty$, then $(W_n)$ is $\Lambda(W)$-upper $\Psi$-regular.
\item Theorem~\ref{thm:compact}: $\{W:\Lambda(W)\le C\}$ is compact in $\delta_\square$.
\item Proposition~\ref{prop:comparison} and Example~\ref{ex:beyond}: $L^p$ upper regularity implies $\Psi$ upper regularity, and sparse random graphs with a limit in $L\log L\setminus\bigcup_{p>1}L^p$ are upper $\Psi$-regular but not upper $L^p$ regular for any $p>1$.
\item Theorem~\ref{thm:stein} and Corollary~\ref{cor:stein}: by exchangeable pairs built from switchings, uniformly random graphs with a given degree sequence satisfy an explicit cut-concentration bound around $\varphi_n\otimes\varphi_n$ with no condition on $d_{\max}$; hence, if the profile converges in $L^1$, $M/n\to\infty$ and $M/n^2\to0$, they converge to $\varphi\otimes\varphi$, and they are upper $\Psi$-regular when $\varphi\in L\log L$. Also: cut-closeness to a $\beta$-model whenever $-\log\pi_x(d)=o(M)$ (Theorem~\ref{thm:degconv}), the McKay-range case with an entropy expansion (Theorem~\ref{thm:degsparse}), and an entrywise Chung--Lu form of the $\beta$-model for $d_{\max}^2=o(M)$ (Proposition~\ref{prop:chunglu}).
\item Proposition~\ref{prop:rate}: for $W=p(1+U)$ with values in $[0,1]$, the rate $I_p(W)$ of $G(n,p)$ satisfies $\frac p2\Lambda(U)\le I_p(W)\le(1+\kappa(p))\frac p2\Lambda(U)$ with $\kappa(p)=O(1/\log\frac1p)$.
\end{itemize}

\subsection{Relation to \cite{BCCZ19}}
We follow the architecture of \cite{BCCZ19}; all proofs are written out in the Orlicz setting. By Lemma~\ref{lem:tail}, upper $\Psi$-regularity implies a uniform tail condition on steppings, and by the de la Vall\'ee-Poussin theorem \cite{RR91} uniform integrability is equivalent to boundedness in the modular of some superlinear Young function. Hence the qualitative part of Theorems~\ref{thm:limits} and~\ref{thm:compact} (without the bound $\Lambda(W)\le C$) is closely related to the uniform $L^1$ theory of \cite[Appendix~C]{BCCZ19}, and may be derivable from it. What is new here is the explicit Orlicz framework with its natural Young pair, the regularity lemma with the modular as energy and polynomial dependence of the iteration count on $1/\varepsilon$, the bound $\Lambda(W)\le C$ for limits, the converse Proposition~\ref{prop:converse}, the examples beyond $L^p$, and the identification of $\Psi$ as the sparse large deviation rate.

\subsection{Motivation}
The class $L\log L$ arises naturally for graphs with a prescribed degree sequence. In the sparse regime a uniformly random graph with degrees $d$ resembles the maximum-entropy $\beta$-model \cite{CDS11}, whose normalized graphon is approximately $\varphi\otimes\varphi$ with $\varphi=d/\bar d$; its entropy involves the profile only through $\int\varphi\log\varphi$, and $\Lambda(\varphi\otimes\varphi)=2\int\varphi\log\varphi+O(1)$ (Section~\ref{sec:degrees}). We prove that such graphs converge to $\varphi\otimes\varphi$ for every sparse $L^1$ profile, with no cap on the maximum degree (Corollary~\ref{cor:stein}), using Stein's method with switchings; we also give a transfer inequality from the $\beta$-model in the spirit of \cite{CDS11,Meh14}. Rank-one graphons with power-law profiles $\mu\{\varphi>t\}\asymp t^{-(\tau-1)}$ lie in some $L^p$, $p>1$, exactly when $\tau>2$; at $\tau=2$ with logarithmic corrections they can lie in $L\log L$ and in no $L^p$ with $p>1$.

\section{Preliminaries}\label{sec:prelim}

\subsection{Graphons, steppings, cut metric}
Let $\mu$ denote Lebesgue measure on $[0,1]$ and on $[0,1]^2$. A \emph{graphon} is a symmetric measurable $W:[0,1]^2\to[0,\infty)$ with $\|W\|_1<\infty$. For a finite measurable partition $\cP=\{J_1,\dots,J_m\}$ of $[0,1]$ into sets of positive measure, the \emph{stepping} $W_\cP$ is the average of $W$ on each $J_a\times J_b$. It is the conditional expectation of $W$ given the $\sigma$-algebra generated by the rectangles $J_a\times J_b$; if $\cP'$ refines $\cP$ then $(W_{\cP'})_\cP=W_\cP$. The \emph{cut norm} and \emph{cut metric} are
\[
\cut{U}=\sup_{S,T\subseteq[0,1]}\Big|\int_{S\times T}U\Big|,\qquad
\delta_\square(U,W)=\inf_\sigma\cut{U^\sigma-W},
\]
where $\sigma$ ranges over measure-preserving bijections of $[0,1]$ and $U^\sigma(x,y)=U(\sigma x,\sigma y)$ \cite{Lov12}.

A simple graph $G$ on $[n]$ defines the graphon $W^G$, equal to $1$ on $I_i\times I_j$ when $ij\in E(G)$ and $0$ otherwise, where $I_i=[\frac{i-1}n,\frac in)$. Its density is $\|W^G\|_1=2|E(G)|/n^2$, and its \emph{normalized graphon} is $W^G/\|W^G\|_1$. A sequence of graphs \emph{converges} to $W$ if $\delta_\square(W^{G_n}/\|W^{G_n}\|_1,W)\to0$.

\subsection{The Young pair}
The function $\Psi(u)=(1+u)\log(1+u)-u$ is convex and increasing on $[0,\infty)$, with $\Psi(0)=0$, $\Psi'(u)=\log(1+u)$, $\Psi''(u)=1/(1+u)$, $\Psi(u)\le u^2/2$, and $\Psi(u)=u\log u-u+O(\log u)$ as $u\to\infty$. Its Legendre conjugate is $\Psi^*(v)=e^v-1-v$, and \eqref{eq:young} is Young's inequality. We write $\Lambda(U)=\int\Psi(|U|)\,d\mu$ and call $W$ an \emph{$L\log L$ graphon} if $\Lambda(W)<\infty$.

\begin{lemma}[Jensen]\label{lem:jensen}
For every graphon $W$ and partition $\cP$, $\Lambda(W_\cP)\le\Lambda(W)$. If $\cP'$ refines $\cP$ then $\Lambda(W_\cP)\le\Lambda(W_{\cP'})$.
\end{lemma}
\begin{proof}
$W_\cP$ is a conditional expectation of $W$ and of $W_{\cP'}$, and $\Psi$ is convex.
\end{proof}

\begin{lemma}[Small sets]\label{lem:smallset}
Let $U\ge0$ be measurable on a probability space and $E$ an event with $\mu(E)<1$. Then
\[
\int_EU\;\le\;\frac{1+\Lambda(U)}{\log\big(1/\mu(E)\big)},
\]
where the right side is read as $0$ if $\mu(E)=0$.
\end{lemma}
\begin{proof}
If $\mu(E)=0$ there is nothing to prove. Otherwise apply \eqref{eq:young} with $v=s\mathbf 1_E$ and integrate: $s\int_EU\le\Lambda(U)+\mu(E)\Psi^*(s)$. Take $s=\log(1/\mu(E))$; then $\mu(E)\Psi^*(s)=1-\mu(E)-\mu(E)s\le1$.
\end{proof}

\begin{lemma}[Tails]\label{lem:tail}
For $U\ge0$ and $K>0$, $\int_{\{U>K\}}U\le\Lambda(U)\,K/\Psi(K)$, and $\Psi(K)/K\ge\log(1+K)-1$.
\end{lemma}
\begin{proof}
$\Psi(u)/u$ is nondecreasing since $\Psi$ is convex with $\Psi(0)=0$, so $u\le\Psi(u)K/\Psi(K)$ for $u>K$. Also $\Psi(K)/K=(1+\frac1K)\log(1+K)-1\ge\log(1+K)-1$.
\end{proof}

\begin{lemma}[Bregman divergence of $\Psi$]\label{lem:bregman}
For $a,b\ge0$ let $D_\Psi(a,b)=\Psi(a)-\Psi(b)-\Psi'(b)(a-b)=(1+a)\log\frac{1+a}{1+b}-(a-b)$. Then
\[
D_\Psi(a,b)\;\ge\;\frac{(a-b)^2}{2(1+\max(a,b))}.
\]
\end{lemma}
\begin{proof}
For fixed $b$, $f(a)=D_\Psi(a,b)$ satisfies $f(b)=f'(b)=0$ and $f''(t)=1/(1+t)\ge1/(1+\max(a,b))$ for $t$ between $a$ and $b$. Integrate twice.
\end{proof}

\subsection{Upper $\Psi$-regularity}

\begin{definition}\label{def:upper}
Let $C\ge0$ and $\eta\in(0,1)$. A graphon $W$ is \emph{$(C,\eta)$-upper $\Psi$-regular} if $\Lambda(W_\cP)\le C$ for every partition $\cP$ of $[0,1]$ into measurable parts of measure at least $\eta$. A simple graph is $(C,\eta)$-upper $\Psi$-regular if its normalized graphon is. A sequence $(W_n)$ (or $(G_n)$) is a \emph{$C$-upper $\Psi$-regular sequence} if for every $\eta>0$ there is $n_0$ such that $W_n$ is $(C+\eta,\eta)$-upper $\Psi$-regular for all $n\ge n_0$.
\end{definition}

The property is invariant under measure-preserving rearrangements, and by Lemma~\ref{lem:jensen} every graphon with $\Lambda(W)\le C$ is $(C,\eta)$-upper $\Psi$-regular for every $\eta$. For a sparse graph the normalized graphon itself has $\Lambda\approx\log(1/\rho)\to\infty$, so upper regularity constrains only coarse steppings, as in \cite{BCCZ19}.

\begin{remark}\label{rem:vertex}
\cite[Def.~2.1]{BCCZ19} defines upper regularity of a graph through partitions of its vertex set, which is a priori weaker than the graphon definition used here. All our results are stated for the graphon definition.
\end{remark}

\begin{lemma}[Total mass]\label{lem:mass}
If $W$ is $(C,\eta)$-upper $\Psi$-regular then $\|W\|_1\le A(C):=\Psi^{-1}(C)$.
\end{lemma}
\begin{proof}
For the trivial partition, $\Psi(\|W\|_1)=\Lambda(W_{\{[0,1]\}})\le C$.
\end{proof}

\begin{proposition}[Comparison with $L^p$]\label{prop:comparison}
Let $1<p\le2$. If $\|W_\cP\|_p\le C$ for all partitions $\cP$ with parts of measure $\ge\eta$, then $W$ is $(C^p/(p-1),\eta)$-upper $\Psi$-regular. For $p>2$ the same hypothesis gives $(C^2/2,\eta)$-upper $\Psi$-regularity.
\end{proposition}
\begin{proof}
$\Psi(u)\le u\log(1+u)$ since $\log(1+u)\le u$. For $q=p-1\in(0,1]$, $\log(1+u)\le((1+u)^q-1)/q\le u^q/q$, so $\Psi(u)\le u^p/(p-1)$. For $p>2$ use $\Psi(u)\le u^2/2$ and $\|\cdot\|_2\le\|\cdot\|_p$.
\end{proof}

\section{A weak regularity lemma}\label{sec:regularity}

\begin{lemma}[Strips]\label{lem:strip}
Let $W$ be $(C,\eta)$-upper $\Psi$-regular, $\cP$ a partition with parts of measure $\ge\eta$, and $X\subseteq[0,1]$ with $\mu(X)\le\eta_0\le\frac12$, where $\eta\le\eta_0$. Then
\[
\int_{X\times[0,1]}W+\int_{X\times[0,1]}W_\cP\;\le\;\frac{2(C+1)}{\log(1/\eta_0)}.
\]
\end{lemma}
\begin{proof}
Choose $X'\supseteq X$ with $\mu(X')=\max(\mu(X),\eta)\le\eta_0$. The partition $\cQ=\{X',[0,1]\setminus X'\}$ has parts of measure $\ge\eta$, and $X'\times[0,1]$ is measurable for $\cQ\times\cQ$, so $\int_{X\times[0,1]}W\le\int_{X'\times[0,1]}W=\int_{X'\times[0,1]}W_\cQ$, using $W\ge0$. Apply Lemma~\ref{lem:smallset} to $W_\cQ$ and to $W_\cP$, both of which have $\Lambda\le C$.
\end{proof}

\begin{theorem}[Weak regularity for upper $\Psi$-regular graphons]\label{thm:regularity}
Let $C>0$ and $0<\varepsilon<1$. Put
\[
A=\Psi^{-1}(C),\qquad N=\Big\lfloor\frac{8C(1+2A)}{\varepsilon^2}\Big\rfloor+1,\qquad \eta_0=e^{-8(C+1)/\varepsilon}.
\]
Let $\cP_0$ be a partition into $m_0$ parts, put $\eta_0'=\min\big(\min_{J\in\cP_0}\mu(J),\ \eta_0/m_0\big)$, and let $W$ be $(C,\eta)$-upper $\Psi$-regular with $\eta\le\eta_0'8^{-N}$. Then there is a refinement $\cP$ of $\cP_0$ into at most $m_04^{N}$ parts, each of measure at least $\eta_0'8^{-N}$, such that
\[
\cut{W-W_\cP}\le\varepsilon\qquad\text{and}\qquad\Lambda(W_\cP)\le C.
\]
\end{theorem}

\begin{proof}
We build refinements $\cP_0\preceq\cP_1\preceq\cdots$ in which $\cP_k$ has $m_k\le m_04^k$ parts, each of measure at least $\eta_k:=\eta_0'8^{-k}$. For $k\le N$ every $\cP_k$ has parts of measure $\ge\eta$, so $0\le\Lambda(W_{\cP_k})\le C$.

Suppose $k<N$ and $\cut{W-W_{\cP_k}}>\varepsilon$; pick $S,T$ with $|\int_{S\times T}(W-W_{\cP_k})|>\varepsilon$.

\emph{Refining without small parts.} Each part $J\in\cP_k$ splits into at most four atoms $J\cap S\cap T,\dots,J\setminus(S\cup T)$, the largest of measure $\ge\mu(J)/4$. Merge every atom of measure $<\eta_k/8$ into the largest atom of the same $J$. The result $\cP_{k+1}$ refines $\cP_k$, has at most $4m_k$ parts, and every part has measure $\ge\eta_k/8=\eta_{k+1}$. Let $S'$ (resp.\ $T'$) be the union of the parts of $\cP_{k+1}$ whose largest constituent atom lies in $S$ (resp.\ $T$). Then $S',T'$ are unions of parts of $\cP_{k+1}$, and $S\triangle S'$, $T\triangle T'$ are contained in the union $X$ of the merged atoms, with
\[
\mu(X)\le3m_k\cdot\frac{\eta_k}8\le m_k\eta_k\le m_0\eta_0'2^{-k}\le\eta_0.
\]

\emph{Discrepancy survives.} Since $(S\times T)\triangle(S'\times T')\subseteq(X\times[0,1])\cup([0,1]\times X)$ and $W,W_{\cP_k}$ are symmetric and nonnegative, Lemma~\ref{lem:strip} gives
\[
\Big|\int_{S\times T}(W-W_{\cP_k})-\int_{S'\times T'}(W-W_{\cP_k})\Big|\le\frac{4(C+1)}{\log(1/\eta_0)}=\frac\varepsilon2.
\]
Since $\mathbf 1_{S'\times T'}$ is $\cP_{k+1}\times\cP_{k+1}$-measurable, $\int_{S'\times T'}W=\int_{S'\times T'}W_{\cP_{k+1}}$, so
\[
\frac\varepsilon2\le\big|\langle W_{\cP_{k+1}}-W_{\cP_k},\mathbf 1_{S'\times T'}\rangle\big|\le\int|W_{\cP_{k+1}}-W_{\cP_k}|.
\]

\emph{Energy increment.} Since $\Psi'(W_{\cP_k})$ is $\cP_k\times\cP_k$-measurable and bounded, and $W_{\cP_k}$ is the conditional expectation of $W_{\cP_{k+1}}$, we have $\int\Psi'(W_{\cP_k})(W_{\cP_{k+1}}-W_{\cP_k})=0$, hence
\[
\Lambda(W_{\cP_{k+1}})-\Lambda(W_{\cP_k})=\int D_\Psi(W_{\cP_{k+1}},W_{\cP_k}).
\]
Write $a=W_{\cP_{k+1}}$, $b=W_{\cP_k}$. By Cauchy--Schwarz with weight $1+\max(a,b)\le1+a+b$ and Lemma~\ref{lem:bregman},
\[
\Big(\int|a-b|\Big)^2\le\int\frac{(a-b)^2}{1+\max(a,b)}\cdot\int(1+a+b)\le2\int D_\Psi(a,b)\cdot(1+2\|W\|_1).
\]
With Lemma~\ref{lem:mass} this gives $\Lambda(W_{\cP_{k+1}})-\Lambda(W_{\cP_k})\ge\varepsilon^2/(8(1+2A))$.

\emph{Termination.} If $\cut{W-W_{\cP_k}}>\varepsilon$ for $k=0,\dots,N-1$, then $\cP_N$ is defined, has parts of measure $\ge\eta$, and $\Lambda(W_{\cP_N})\ge N\varepsilon^2/(8(1+2A))>C$, a contradiction. Hence $\cut{W-W_{\cP_k}}\le\varepsilon$ for some $k<N$, and $\cP=\cP_k$ has the stated properties.
\end{proof}

\begin{remark}
The number of parts is $4^{O(C(1+A)/\varepsilon^2)}$, as in the classical weak regularity lemma; the bound on $\eta$ is $e^{-O(C/\varepsilon)}8^{-N}$, the factor $e^{-O(C/\varepsilon)}$ coming from the logarithmic decay in Lemma~\ref{lem:smallset}. Using the modular itself as the energy avoids the truncation used for $1<p<2$ in \cite[Thm.~4.3]{BCCZ19}.
\end{remark}

\begin{remark}[General Young functions]\label{rem:generalPsi}
Let $\Phi$ be a Young function with $\Phi(u)/u\to\infty$ and conjugate $\Phi^*$. Young's inequality gives $\int_EU\le\omega_\Phi(\mu(E);\Lambda_\Phi(U))$ with $\omega_\Phi(t;L)=\inf_{s>0}(L+t\Phi^*(s))/s\to0$ as $t\to0$. Lemma~\ref{lem:strip}, and hence the discrepancy transfer, holds with this modulus. The energy step uses $\Phi''$ and may be replaced by the truncated-energy argument of \cite[Thm.~4.3]{BCCZ19} with Lemma~\ref{lem:tail} in place of H\"older. Superlinearity of $\Phi$ is used in these two estimates and in the uniform integrability step of Theorem~\ref{thm:limits}; for $\Phi(u)=u$ all three fail.
\end{remark}

\section{Limits of upper $\Psi$-regular sequences}\label{sec:limits}

\begin{theorem}[Existence of limits]\label{thm:limits}
Let $(W_n)$ be a $C$-upper $\Psi$-regular sequence of graphons. Then there is a graphon $W$ with $\Lambda(W)\le C$ and a subsequence with $\delta_\square(W_{n_j},W)\to0$. The same holds for a $C$-upper $\Psi$-regular sequence of simple graphs, with $W_n$ their normalized graphons.
\end{theorem}

\begin{proof}
\emph{Nested partitions.} Let $\varepsilon_k=1/(k+1)$. Apply Theorem~\ref{thm:regularity} with constant $C+1$ repeatedly: first with $\cP_0=\{[0,1]\}$ and $\varepsilon_1$, then starting from the partition just obtained with $\varepsilon_2$, and so on. At stage $k$ the number of parts $M_k$ and the lower bound $\nu_k>0$ on their measures depend only on $C$ and $k$, and the theorem applies once $W_n$ is $(C+1,\eta_k)$-upper $\Psi$-regular for some $\eta_k=\eta_k(C)>0$, which holds for $n\ge n_0(k)$. For such $n$ we obtain nested partitions $\cP_{n,1}\preceq\cdots\preceq\cP_{n,k}$ with $|\cP_{n,k}|\le M_k$, parts of measure $\ge\nu_k$, and $\cut{W_n-(W_n)_{\cP_{n,k}}}\le\varepsilon_k$. Moreover, for every $\eta'\in(0,\nu_k]$, $\Lambda((W_n)_{\cP_{n,k}})\le C+\eta'$ for all large $n$, so $\limsup_n\Lambda((W_n)_{\cP_{n,k}})\le C$.

\emph{Interval representatives.} Order the parts of $\cP_{n,1}$, then within each part order its children in $\cP_{n,2}$, and so on, and lay them out as consecutive intervals of the same measures. This gives nested interval partitions $\cI_{n,k}$ and step functions $U_{n,k}$ on $\cI_{n,k}\times\cI_{n,k}$ with the same block values as $(W_n)_{\cP_{n,k}}$. Since $(W_n)_{\cP_{n,k}}$ and $U_{n,k}$ are step functions with the same block structure and block measures, $\delta_\square((W_n)_{\cP_{n,k}},U_{n,k})=0$, and hence $\delta_\square(W_n,U_{n,k})\le\varepsilon_k$. Also $U_{n,k}=\E[U_{n,k+1}\mid\cI_{n,k}\times\cI_{n,k}]$.

\emph{Compactness.} For fixed $k$ the data of $U_{n,k}$ are the number of parts ($\le M_k$), at most $M_k$ interval lengths in $[\nu_k,1]$, and at most $M_k^2$ block values in $[0,A(C+1)/\nu_k^2]$ by Lemma~\ref{lem:mass}. By a diagonal argument we pass to a subsequence along which, for every $k$, the number of parts is constant and all lengths and values converge. Let $\cI_k$ be the limiting interval partition (its intervals have length $\ge\nu_k$) and $U_k$ the limiting step function. Then $U_{n,k}\to U_k$ in $L^1$, the partitions $\cI_k$ are nested, $U_k=\E[U_{k+1}\mid\cI_k\times\cI_k]$, and, since $\Lambda$ of a step function is continuous in its finitely many parameters, $\Lambda(U_k)\le C$.

\emph{Martingale.} $(U_k)$ is a martingale for the filtration generated by $\cI_k\times\cI_k$, bounded in the modular of the superlinear Young function $\Psi$. By the de la Vall\'ee-Poussin theorem it is uniformly integrable, so it converges almost surely and in $L^1$ to some $U$ with $U_k=\E[U\mid\cI_k\times\cI_k]$. By Fatou, $\Lambda(U)\le C$; $U$ is symmetric and nonnegative, so it is a graphon.

\emph{Convergence.} For each $k$ and $n$ in the subsequence,
\[
\delta_\square(W_n,U)\le\delta_\square(W_n,U_{n,k})+\|U_{n,k}-U_k\|_1+\|U_k-U\|_1\le\varepsilon_k+o_n(1)+\|U_k-U\|_1.
\]
Letting $n\to\infty$ and then $k\to\infty$ gives $\delta_\square(W_n,U)\to0$. The statement for graphs is the special case of normalized graphons.
\end{proof}

\begin{remark}
For $\Phi(u)=u$ the conclusion fails: \cite[Prop.~2.12]{BCCZ19} gives normalized simple graph sequences that are $\delta_\square$-Cauchy without a graphon limit.
\end{remark}

\begin{proposition}[Convergent sequences are upper regular]\label{prop:converse}
Let $W$ be a graphon with $\Lambda(W)<\infty$ and $(W_n)$ graphons with $\delta_\square(W_n,W)\to0$. Then $(W_n)$ is a $\Lambda(W)$-upper $\Psi$-regular sequence.
\end{proposition}
\begin{proof}
Fix $\eta>0$ and choose $\sigma_n$ with $t_n:=\cut{W_n-W^{\sigma_n}}\to0$. Let $\cP$ have parts of measure $\ge\eta$. On each block $J_a\times J_b$, $|(W_n)_\cP-(W^{\sigma_n})_\cP|\le t_n/(\mu(J_a)\mu(J_b))\le t_n/\eta^2$. Since $\Psi$ is increasing and convex, $\Psi(u+t)\le\Psi(u)+t\log(1+u+t)$, and $(W^{\sigma_n})_\cP\le\|W\|_1/\eta^2$, so
\[
\Lambda((W_n)_\cP)\le\Lambda((W^{\sigma_n})_\cP)+\frac{t_n}{\eta^2}\log\Big(1+\frac{\|W\|_1+t_n}{\eta^2}\Big)\le\Lambda(W)+\frac{t_n}{\eta^2}\log\Big(1+\frac{\|W\|_1+t_n}{\eta^2}\Big),
\]
by Lemma~\ref{lem:jensen} and invariance of $\Lambda$ under $\sigma_n$. The error is uniform in $\cP$ and tends to $0$ for fixed $\eta$.
\end{proof}

\begin{theorem}[Compactness]\label{thm:compact}
For every $C\ge0$ the set $\{W:\Lambda(W)\le C\}$ of graphons, with graphons at cut distance $0$ identified, is compact in $\delta_\square$.
\end{theorem}
\begin{proof}
A sequence in this set is $C$-upper $\Psi$-regular by Lemma~\ref{lem:jensen}; Theorem~\ref{thm:limits} gives a subsequence converging to a graphon with $\Lambda\le C$.
\end{proof}

\section{Sparse random graphs and examples beyond $L^p$}\label{sec:examples}

For a graphon $W$, $n\ge1$ and $\rho>0$, let $\mathbf G(n,W,\rho)$ be the random graph on $[n]$ obtained by sampling $x_1,\dots,x_n$ uniformly and joining $i\ne j$ independently with probability $\min(1,\rho W(x_i,x_j))$. By \cite[Thm.~2.14, Cor.~2.15]{BCCZ19}, for every $L^1$ graphon with $\|W\|_1>0$, if $\rho_n\to0$ and $n\rho_n\to\infty$, then almost surely $\mathbf G(n,W,\rho_n)$ converges to $W/\|W\|_1$ after normalization. With Proposition~\ref{prop:converse}:

\begin{corollary}\label{cor:wrandom}
If $\Lambda(W/\|W\|_1)<\infty$, $\rho_n\to0$ and $n\rho_n\to\infty$, then almost surely $(\mathbf G(n,W,\rho_n))_n$ is a $\Lambda(W/\|W\|_1)$-upper $\Psi$-regular sequence.
\end{corollary}

\begin{example}[Beyond every $L^p$]\label{ex:beyond}
Let $\varphi(x)=2x^{-1}(\log\frac ex)^{-3}$ on $(0,1]$; the substitution $t=\log\frac ex$ shows $\int_0^1\varphi=1$.
\begin{enumerate}[leftmargin=2em,label=(\alph*)]
\item $\varphi\log\varphi\asymp x^{-1}(\log\frac ex)^{-2}$ near $0$, so $\varphi\in L\log L$; but $\varphi\notin L^p$ for any $p>1$.
\item $W=\varphi\otimes\varphi$ satisfies $\Lambda(W)\le2\int\varphi\log(1+\varphi)<\infty$, since $\Psi(uv)\le uv\log(1+uv)\le uv(\log(1+u)+\log(1+v))$; and $W\notin L^p$ for any $p>1$.
\item By Corollary~\ref{cor:wrandom}, $G_n=\mathbf G(n,W,\rho_n)$ with $\rho_n\to0$ and $n\rho_n\to\infty$ is almost surely a $\Lambda(W)$-upper $\Psi$-regular sequence converging to $W$.
\item Almost surely, $(G_n)$ is not a $C$-upper $L^p$ regular sequence (in the graphon sense: $\|(W_n)_\cP\|_p\le C+\eta$ for parts of measure $\ge\eta$ and $n\ge n_0(\eta)$) for any $p>1$ and $C<\infty$. Indeed, fix $\eta$ and a partition $\cP$ with parts of measure $\ge\eta$, and let $\sigma_n$ with $t_n=\cut{W_n-W^{\sigma_n}}\to0$. Applying the hypothesis to $\sigma_n^{-1}\cP$ and the blockwise bound $|(W_n)_{\sigma_n^{-1}\cP}-(W^{\sigma_n})_{\sigma_n^{-1}\cP}|\le t_n/\eta^2$ gives $\|W_\cP\|_p\le C+\eta+t_n/\eta^2$ for all large $n$, hence $\|W_\cP\|_p\le C+\eta$. For the interval partitions $\cP_k$ into $2^k$ equal parts (take $\eta=2^{-k}$), $\|W_{\cP_k}\|_p\le C+2^{-k}$ and $W_{\cP_k}\to W$ almost surely by martingale convergence, so Fatou gives $\|W\|_p\le C$, a contradiction.
\end{enumerate}
In degree terms, $\mu\{\varphi>t\}\asymp t^{-1}(\log t)^{-3}$: a power law with exponent $\tau=2$ and a logarithmic correction.
\end{example}

\begin{remark}
For the vertex-partition definition of \cite[Def.~2.1]{BCCZ19}, part (d) should follow in the same way by grouping vertices according to the part of $\cP$ containing $x_i$ and using concentration of the block edge counts, whose expectations tend to infinity; we do not write this out.
\end{remark}

\begin{remark}[Power laws]
For connection probabilities $\min(1,n^\beta(ij)^{-\alpha})$ with $0<\alpha<1$ \cite{BCCZ19}, the limit $(xy)^{-\alpha}$ lies in $L^p$ for $1<p<1/\alpha$, so $L\log L$ adds nothing there. The new range is the endpoint $\alpha\to1$, i.e.\ degree exponent $\tau\to2$, with logarithmic corrections as in Example~\ref{ex:beyond}.
\end{remark}

\section{Graphs with a given degree sequence}\label{sec:degrees}

Throughout this section $d=d^{(n)}$ is a graphical degree sequence on $[n]$ with $M=\frac12\sum_id_i$ edges, $d_{\max}=\max_id_i$, density $\rho=2M/n^2$, average degree $\bar d=2M/n$, and normalized profile $\varphi_n(x)=d_{\lceil nx\rceil}/\bar d$, so that $\int\varphi_n=1$. Let $G_d$ be a uniformly random graph with degree sequence $d$.

\subsection{$\beta$-models and the transfer inequality}
For $x\in(0,\infty)^n$ the \emph{$\beta$-model} $G_x$ has independent edges with probabilities $w_{ij}=x_ix_j/(1+x_ix_j)$. Every graph with degree sequence $d$ has the same probability $\prod_{i<j}(1-w_{ij})\prod_ix_i^{d_i}$ under $G_x$ \cite{CDS11,Meh14}, so $G_x$ conditioned on $\deg(G_x)=d$ is uniform. Writing $\pi_x(d)=\Pr[\deg(G_x)=d]$, it follows that for every event $\mathcal A$,
\begin{equation}\label{eq:transfer}
\Pr[G_d\in\mathcal A]\;\le\;\Pr[G_x\in\mathcal A]/\pi_x(d).
\end{equation}
The maximum-entropy choice $x^\star$, for which $\E\deg(G_{x^\star})=d$, maximizes $\pi_x(d)$ (for strictly graphical $d$, $\pi_x(d)=\pi_{x^\star}(d)\,e^{-\sum_{i<j}\mathrm{KL}(w^\star_{ij}\|w_{ij})}$); the conjecture of \cite[Conj.~2.18]{Meh14} asserts that $-\log\pi_{x^\star}(d)=O(n\log n)$ for strictly graphical $d$.

\subsection{Cut-norm concentration}
For a symmetric array $w=(w_{ij})$ with zero diagonal let $W^w$ be the step graphon with value $w_{ij}$ on $I_i\times I_j$.

\begin{lemma}\label{lem:cutconc}
Let $G$ have independent edges with probabilities $w_{ij}$, $\sum_{i<j}w_{ij}\le M$. For $0<\varepsilon\le1$,
\[
\Pr\big[\cut{W^G-W^w}\ge\varepsilon\rho\big]\;\le\;2\cdot4^n\,e^{-\varepsilon^2M/3}.
\]
\end{lemma}
\begin{proof}
$W^G-W^w$ is a step function on $\{I_i\}$, and $\int_{S\times T}$ of such a function is multilinear in the fractions $\mu(S\cap I_i)n,\mu(T\cap I_j)n\in[0,1]$, so the supremum defining the cut norm is attained when $S,T$ are unions of intervals $I_i$. For vertex sets $S,T\subseteq[n]$, $n^2\int_{S\times T}(W^G-W^w)=X:=\sum_{i<j}c_{ij}(A_{ij}-w_{ij})$ with $c_{ij}=\mathbf 1[i\in S,j\in T]+\mathbf 1[j\in S,i\in T]\in\{0,1,2\}$. The summands are independent, bounded by $2$, and $\sum_{i<j}\Var(c_{ij}A_{ij})\le4M$. Bernstein's inequality with $t=\varepsilon\rho n^2=2\varepsilon M$ gives $\Pr[|X|\ge t]\le2\exp(-t^2/(2(4M+2t/3)))\le2e^{-\varepsilon^2M/3}$. Take a union bound over the $4^n$ pairs $(S,T)$.
\end{proof}

\subsection{Convergence}

\begin{theorem}[Transfer of convergence]\label{thm:degconv}
Let $x=x^{(n)}$ be $\beta$-model parameters with $\sum_{i<j}w_{ij}\le M$, and suppose $M/n\to\infty$ and $-\log\pi_x(d)=o(M)$. Then
\[
\cut{W^{G_d}/\rho-W^w/\rho}\;\to\;0\qquad\text{in probability}.
\]
\end{theorem}
\begin{proof}
By \eqref{eq:transfer} and Lemma~\ref{lem:cutconc}, for fixed $\varepsilon\in(0,1]$ the probability that the cut norm is at least $\varepsilon$ is at most $2\exp(n\log4-\varepsilon^2M/3-\log\pi_x(d))\to0$.
\end{proof}

With the maximum-entropy parameters, Theorem~\ref{thm:degconv} says: whenever $-\log\pi_{x^\star}(d)=o(M)$ and $M/n\to\infty$, a uniformly random graph with degree sequence $d$ is cut-close, after normalization, to its maximum-entropy $\beta$-model. In particular \cite[Conj.~2.18]{Meh14} implies this for all strictly graphical sequences with $M/(n\log n)\to\infty$. In the dense regime this is the mechanism of \cite{CDS11}, where a bound of the form $\pi(d)\ge e^{-o(n^2)}$ is used. In the sparse regime we can take the explicit Chung--Lu parameters and obtain an unconditional result.

\begin{theorem}[Sparse degree-constrained graphs converge]\label{thm:degsparse}
Suppose $M\to\infty$, $M/n\to\infty$ and $d_{\max}=o(M^{1/4})$. Then
\[
\cut{W^{G_d}/\rho-\varphi_n\otimes\varphi_n}\to0\qquad\text{in probability}.
\]
Consequently, if $\varphi_n\to\varphi$ in $L^1$, then $G_d\to\varphi\otimes\varphi$ in probability in the normalized cut metric; every subsequence has a further subsequence along which the convergence is almost sure, and along it $(G_d)$ is a $\Lambda(\varphi\otimes\varphi)$-upper $\Psi$-regular sequence whenever $\Lambda(\varphi\otimes\varphi)<\infty$ (Proposition~\ref{prop:converse}).
\end{theorem}

\begin{proof}
Take $x_i=d_i/\sqrt{2M}$, so $w_{ij}=u_{ij}/(1+u_{ij})$ with $u_{ij}=d_id_j/(2M)\le d_{\max}^2/(2M)\to0$, and $\sum_{i<j}w_{ij}\le\sum_{i<j}u_{ij}\le M$.

\emph{The $\beta$-model is close to $\varphi_n\otimes\varphi_n$.} On $I_i\times I_j$ with $i\ne j$, $\varphi_n\otimes\varphi_n=u_{ij}/\rho$ and $0\le(u_{ij}-w_{ij})/\rho\le u_{ij}^2/\rho$; on the diagonal blocks $W^w=0$. Hence
\[
\|\varphi_n\otimes\varphi_n-W^w/\rho\|_1\le\frac1{n^2\rho}\Big(\sum_{i\ne j}u_{ij}^2+\sum_iu_{ii}\Big)\le\frac{d_{\max}^2}{2M}+\frac{d_{\max}}{2M}\to0,
\]
using $\sum_{i,j}u_{ij}^2\le(\sum_id_i^2)^2/(4M^2)\le d_{\max}^2$ and $\sum_iu_{ii}\le d_{\max}$.

\emph{Lower bound on $\pi_x(d)$.} We have $\log\pi_x(d)=\log G(d)-\sum_{i<j}\log(1+u_{ij})+\sum_id_i\log x_i$, where $G(d)$ is the number of graphs with degree sequence $d$. By McKay's formula \cite{McK85}, valid for $d_{\max}=o(M^{1/4})$,
\[
G(d)=\frac{(2M)!}{M!\,2^M\prod_id_i!}\exp\big(-\lambda-\lambda^2+o(1)\big),\qquad\lambda=\frac1{4M}\sum_id_i(d_i-1)\le\frac{d_{\max}}2.
\]
Using $\log(1+u)\le u$, $\sum_{i<j}u_{ij}=M-\sum_id_i^2/(4M)$, the Stirling bounds $\log\frac{(2M)!}{M!2^M}\ge M\log(2M)-M-\log M-1$ and $\log k!\le k\log k-k+\log k+1$ $(k\ge1)$, and Jensen's inequality $\sum_i\log d_i\le n\log\bar d$, the terms $\pm\sum_id_i\log d_i$ and $\pm M\log(2M)$ cancel and
\[
-\log\pi_x(d)\le n\big(1+\log\bar d\big)+\log M+\lambda^2+2+\Big(\lambda-\frac{\sum_id_i^2}{4M}\Big)\le n(1+\log\bar d)+\log M+\frac{d_{\max}^2}4+2
\]
for large $n$. Since $\bar d=2M/n\to\infty$, $n\log\bar d=o(M)$; also $d_{\max}^2=o(M^{1/2})$. So $-\log\pi_x(d)=o(M)$.

Theorem~\ref{thm:degconv} and the triangle inequality give the first claim. The remaining claims follow since $\|\varphi_n\otimes\varphi_n-\varphi\otimes\varphi\|_1\le2\|\varphi_n-\varphi\|_1$ and convergence in probability implies almost sure convergence along a subsequence of any subsequence.
\end{proof}

\begin{remark}
The range $d_{\max}=o(M^{1/4})$ comes only from McKay's formula. Theorem~\ref{thm:stein} below removes every restriction on $d_{\max}$ for cut convergence, by a different route. Local statistics such as tree densities were compared with the $\beta$-model in \cite{Meh14}, unconditionally for relaxed degree constraints $|\deg_i-d_i|<d_i^a$, $a>\frac12$, and conditionally on \cite[Conj.~2.18]{Meh14} for exact degrees.
\end{remark}

\subsection{The sparse $\beta$-model is the Chung--Lu model}\label{subsec:chunglu}

\begin{proposition}[Entrywise Chung--Lu approximation]\label{prop:chunglu}
Let $d$ be strictly graphical with all $d_i\ge1$, $M\ge d_{\max}$, and $\varepsilon:=d_{\max}^2/M\le1/40$. Let $x$ be the maximum-entropy parameters, $y_i=d_i/\sqrt{2M}$, and $q_{ij}=y_iy_j/(1+y_iy_j)=d_id_j/(2M+d_id_j)$. Then for all $i$ and all $i\ne j$,
\[
1-10\varepsilon\le\frac{x_i}{y_i}\le1+19\varepsilon,\qquad\Big|\frac{w_{ij}}{q_{ij}}-1\Big|\le48\varepsilon .
\]
\end{proposition}
\begin{proof}
Write $f(t)=t/(1+t)$, $S=\sum_jx_j$, $m=\max_jx_j=x_k$.
(i) $2M=\sum_{i\ne j}f(x_ix_j)\le S^2$, and $d_j\le x_jS$, so $x_j\ge d_j/S$.
(ii) Since $f$ is increasing and concave, with $T=md_{\max}/S\ge md_j/S$,
\[
d_k=\sum_{j\ne k}f(mx_j)\ge\sum_{j\ne k}f\Big(\frac{md_j}S\Big)\ge\frac{f(T)}T\cdot\frac mS(2M-d_k)=\frac{f(T)(2M-d_k)}{d_{\max}},
\]
so $f(T)\le d_{\max}^2/(2M-d_{\max})\le\varepsilon$, hence $T\le2\varepsilon$ and $m\le2\varepsilon S/d_{\max}\le2\varepsilon S$.
(iii) From $d_j\ge x_j(S-x_j)/(1+x_jm)$ we get $x_j\le d_j(1+m^2)/(S-m)$; summing, $S(S-m)\le2M(1+m^2)$. With $m^2\le4\varepsilon^2S^2/d_{\max}^2=4d_{\max}^2S^2/M^2$ and $S-m\ge(1-2\varepsilon)S$ this gives $S^2(1-10\varepsilon)\le2M$, and then $m^2\le8\varepsilon/(1-10\varepsilon)\le16\varepsilon$.
(iv) Hence $x_i\ge d_i/S\ge(1-10\varepsilon)y_i$ and $x_i\le d_i(1+16\varepsilon)/((1-2\varepsilon)S)\le(1+19\varepsilon)y_i$, using $S\ge\sqrt{2M}$. So $x_ix_j/(y_iy_j)\in[1-20\varepsilon,1+48\varepsilon]$, and since $f(t)/f(s)$ lies between $1$ and $t/s$, the same bounds hold for $w_{ij}/q_{ij}$.
\end{proof}

Thus, whenever $d_{\max}^2=o(M)$, the maximum-entropy multipliers are explicit up to $1+o(1)$: $\lambda_i=\log\frac{\sqrt{2M}}{d_i}+O(\varepsilon)$ uniformly in $i$, and the $\beta$-model is entrywise the model $q_{ij}=d_id_j/(2M+d_id_j)$ of \cite[Thm.~2.25]{Meh14}. That the two agree asymptotically is noted in \cite{GO24}; Proposition~\ref{prop:chunglu} gives an explicit error. Beyond $d_{\max}\approx\sqrt M$ hub pairs saturate, and the entrywise statement fails.

\subsection{Concentration by exchangeable pairs: no cap on $d_{\max}$}\label{subsec:stein}

The transfer inequality pays $-\log\pi(d)$. Stein's method of exchangeable pairs \cite{Cha07}, with switchings as the pair as in \cite{Joh15,CGJ18}, works directly on $G_d$ and pays nothing. For $S,T\subseteq[n]$ let $e(S,T)=\sum_{i\in S,j\in T}A_{ij}$ (ordered pairs), $D_S=\sum_{i\in S}d_i$, and for $k\ge0$ let $J_k$ be the sum of the $k$ largest degrees. Put
\[
K(d)=\frac1{2M}\sum_{a}d_aJ_{d_a},\qquad \mu_d(L)=\frac{J_L}{2M},\qquad \tau_d(L)=\frac1{2M}\sum_{d_a>L}d_a .
\]
Thus $\mu_d(L)$ is the share of degree mass on the $L$ largest vertices and $\tau_d(L)$ the share on vertices of degree $>L$; splitting the sum defining $K$ at $d_a=L$ gives $K(d)/M\le2\mu_d(L)+2\tau_d(L)$ for every $L$.

\emph{The swap chain.} Given a graph $G$ with degrees $d$, choose two oriented edges $(a,b),(c,e)$ independently and uniformly among the $2M$ orientations. If $a,b,c,e$ are distinct and $ac,be\notin E(G)$, let $G'=G-ab-ce+ac+be$; otherwise $G'=G$. The map $((a,b),(c,e))\mapsto((a,c),(b,e))$ is an involution between the proposals taking $G$ to $G'$ and those taking $G'$ to $G$, so the kernel is symmetric; with $G\sim G_d$ uniform, $(G,G')$ is exchangeable.

\begin{lemma}[Drift]\label{lem:drift}
Let $h(G)=e(S,T)-D_SD_T/(2M)$, $\lambda=2/M$, and $r(G)$ the rejection probability. Then $E[h(G)-h(G')\mid G]=\lambda h(G)+\varrho(G)$ with $|\varrho|\le8r$, and for every $G$,
\[
r(G)\le\bar r:=\frac{2d_{\max}+K(d)}{M}.
\]
\end{lemma}
\begin{proof}
Ignore rejections first. The tail $a$ and head $b$ of a uniform oriented edge each have law $d_i/(2M)$, and $\Pr[a\in S,b\in T]=e(S,T)/(2M)$. Removing $ab$ and $ce$ changes $e(S,T)$ by $-2e(S,T)/M$ in expectation, and adding $ac$ and $be$, whose endpoints are independent with law $d_i/(2M)$, by $+D_SD_T/M^2$. So the unconstrained drift of $h(G')-h(G)$ is exactly $-\lambda h(G)$. A rejected proposal would have changed $h$ by at most $8$, which gives $|\varrho|\le8r$. Rejection needs a coincidence among $a,b,c,e$, of probability at most $4\cdot d_{\max}/(2M)$, or $ac\in E$ or $be\in E$. For independent degree-biased $a,c$,
\[
\Pr[ac\in E]=\sum_a\frac{d_a}{2M}\sum_{c\sim a}\frac{d_c}{2M}\le\sum_a\frac{d_a}{2M}\cdot\frac{J_{d_a}}{2M}=\frac{K(d)}{2M},
\]
since the $d_a$ neighbours of $a$ have total degree at most $J_{d_a}$, whatever $G$ is. The same holds for $be$.
\end{proof}

\begin{theorem}[Cut concentration without a degree cap]\label{thm:stein}
For every graphical $d$, every $S,T\subseteq[n]$ and $t>0$,
\[
\Pr\Big[\,\big|e(S,T)-\tfrac{D_SD_T}{2M}\big|\ge4M\bar r+t\Big]\le2e^{-t^2/(32M)}.
\]
Consequently, with probability at least $1-2e^{-n}$, for every $L\ge1$,
\[
\big\|W^{G_d}/\rho-\varphi_n\otimes\varphi_n\big\|_\square\le4\big(\mu_d(L)+\tau_d(L)\big)+\frac{5d_{\max}}M+5\sqrt{\frac nM}.
\]
\end{theorem}
\begin{proof}
Let $m(\theta)=Ee^{\theta h}$, $\theta\ge0$. By the drift lemma, $E[he^{\theta h}]=\lambda^{-1}E[(h-h')e^{\theta h}]-\lambda^{-1}E[\varrho e^{\theta h}]$. By exchangeability, $E[(h-h')e^{\theta h}]=\frac12E[(h-h')(e^{\theta h}-e^{\theta h'})]$. Using $(x-y)(e^x-e^y)\le\frac12(x-y)^2(e^x+e^y)$, exchangeability again, and $|h-h'|\le8$, this is at most $32\theta\,m(\theta)$. Hence $m'(\theta)\le(16M\theta+4M\bar r)m(\theta)$, so $\log m(\theta)\le8M\theta^2+4M\bar r\theta$. Markov's inequality with $\theta=t/(16M)$ gives $\Pr[h\ge4M\bar r+t]\le e^{-t^2/(32M)}$. The lower tail is the same with $-h$.

For the cut norm, $\int_{S\times T}W^{G_d}/\rho=e(S,T)/(2M)$ and $\int_{S\times T}\varphi_n\otimes\varphi_n=D_SD_T/(4M^2)$, up to the diagonal blocks, which carry at most $d_{\max}/(2M)$. The cut norm of a step function is attained on vertex sets. Take $t=\sqrt{32M(2n\log2+n)}\le8.8\sqrt{Mn}$ and a union bound over the $4^n$ pairs. Then $4M\bar r/(2M)\le4d_{\max}/M+2K/M\le4d_{\max}/M+4(\mu_d(L)+\tau_d(L))$.
\end{proof}

\begin{corollary}[Every sparse $L^1$ profile]\label{cor:stein}
If $\varphi_n\to\varphi$ in $L^1$, $M/n\to\infty$ and $\rho=2M/n^2\to0$, then $G_d\to\varphi\otimes\varphi$ in the normalized cut metric, almost surely along any sequence, and $(G_d)$ is $\Lambda(\varphi\otimes\varphi)$-upper $\Psi$-regular whenever $\int\varphi\log\varphi<\infty$ (Proposition~\ref{prop:converse}, Lemma~\ref{lem:rankone}). No cap on $d_{\max}$, no minimum degree, no tameness and no enumeration is needed.
\end{corollary}
\begin{proof}
$d_{\max}\le n=o(M)$. Take $L=\lceil\sqrt M\rceil$. The $L$ largest vertices occupy a fraction $L/n\approx\sqrt{\rho/2}\to0$ of $[0,1]$, and degree $>L$ means $\varphi_n>\sqrt M/\bar d=1/\sqrt{2\rho}\to\infty$. Since $(\varphi_n)$ is uniformly integrable, $\mu_d(L),\tau_d(L)\to0$. Borel--Cantelli handles the $2e^{-n}$.
\end{proof}

\begin{remark}[Prior work, and why no cap is needed]\label{rem:steinprior}
Switchings were used as exchangeable pairs by Johnson \cite{Joh15}, and as size-biased couplings by Cook, Goldstein and Johnson \cite{CGJ18}, who prove concentration of $e(S,T)$ for random $d$-regular graphs and deduce the spectral gap; Tikhomirov and Youssef \cite{TY19} treat near-regular sequences. Gao and Ohapkin \cite[Thm.~1.7]{GO24} embed $\mathcal G(n,(1-o(1))P)$, $P_{ij}=d_id_j/(2M+d_id_j)$, in $G_d$ when $d_{\max}^2=o(M)$ and the minimum degree is $\gg\log n$, which gives Corollary~\ref{cor:stein} in that range. Theorem~\ref{thm:stein} is a heterogeneous, cap-free extension of \cite{CGJ18}. The bound $K(d)/M$ on the rejection rate holds for every graph, and it tends to $0$ for every uniformly integrable profile as $\rho\to0$. Hub pairs with $d_id_j>M$ do saturate, but they carry a vanishing share $2\int_{\varphi>c\rho^{-1/2}}\varphi$ of the mass, and a vertex of degree $\Theta(n)$ is harmless once its share $\Theta(n)/M$ vanishes. In the dense regime $K(d)/M\not\to0$; the $\beta$-model is then the natural limit, and it is reversible for the swap chain, since a swap $ab,ce\to ac,be$ multiplies its weight by $x_ax_cx_bx_e/(x_ax_bx_cx_e)=1$. Convergence to the $\beta$-model at every density with $M\gg n\log n\log(M/n)$ is proved in \cite{Meh-deg} by a counting argument and Theorem~\ref{thm:degconv}. Simulations with a tail-exponent-$2.5$ profile and $d_{\max}$ between $7\sqrt M$ and $21\sqrt M$ give cut deviations of about $0.03$--$0.05$.
\end{remark}

\subsection{Local concentration and regularity relative to the degrees}\label{subsec:relreg}

Theorem~\ref{thm:stein} uses the crude variance bound $|h-h'|\le8$. Using the true conditional variance gives a Bernstein-type bound whose scale is the local edge count, which is what sparse (Kohayakawa--R\"odl) regularity needs. Write $\mu_{S,T}=D_SD_T/(2M)$ and $\Delta_d=4(2d_{\max}+K(d))$.

\begin{theorem}[Local concentration]\label{thm:localstein}
For every graphical $d$, every $S,T\subseteq[n]$ and $t>0$,
\begin{align*}
\Pr\big[e(S,T)-\mu_{S,T}\ge\Delta_d+t\big]&\le\exp\Big(-\frac{t^2}{256(\mu_{S,T}+\Delta_d+t)}\Big),\\
\Pr\big[e(S,T)-\mu_{S,T}\le-\Delta_d-t\big]&\le\exp\Big(-\frac{t^2}{128\,\mu_{S,T}}\Big).
\end{align*}
\end{theorem}
\begin{proof}
Let $h=e(S,T)-\mu_{S,T}$ and $\mu=\mu_{S,T}$. A proposal changes $e(S,T)$ only if a removed edge has an orientation in $S\times T$ (probability at most $e(S,T)/M$ for each of the two edges) or an added pair meets $S\times T\cup T\times S$ (probability at most $D_SD_T/(2M^2)$ for each). Hence, since $|h-h'|\le8$,
\[
E[(h-h')^2\mid G]\le64\Big(\frac{2e(S,T)}M+\frac{D_SD_T}{M^2}\Big)=\frac{128}M(h+2\mu).
\]
For $\theta\ge0$ and $m(\theta)=Ee^{\theta h}$, as in the proof of Theorem~\ref{thm:stein},
\[
m'(\theta)\le\tfrac1\lambda\cdot\tfrac\theta2E\big[(h-h')^2e^{\theta h}\big]+\Delta_dm\le32\theta\big(m'+2\mu m\big)+\Delta_dm .
\]
So $(\log m)'\le(64\theta\mu+\Delta_d)/(1-32\theta)$ for $\theta<1/32$. For $\theta\le1/64$ this gives $\log m(\theta)\le\Delta_d\theta+64\theta^2(\mu+\Delta_d)$. Take $\theta=\min\{1/64,\,t/(128(\mu+\Delta_d))\}$. If $t\le2(\mu+\Delta_d)$ the exponent is at most $-t^2/(256(\mu+\Delta_d))$, and otherwise it is at most $-t/128\le-t^2/(256(\mu+\Delta_d+t))$. For the lower tail apply the same argument to $g=-h$: now $E[(h-h')^2\mid G]\le\frac{128}M(2\mu-g)$, which gives $(1+32\theta)m_g'\le(64\theta\mu+\Delta_d)m_g$ and $\log m_g(\theta)\le32\mu\theta^2+\Delta_d\theta$; take $\theta=t/(64\mu)$.
\end{proof}

\begin{definition}
For $V,V'\subseteq[n]$, the pair $(V,V')$ is \emph{$\varepsilon$-regular relative to $d$} if $|e(S,T)-\mu_{S,T}|\le\varepsilon\,\mu_{V,V'}$ for all $S\subseteq V$ and $T\subseteq V'$.
\end{definition}

This is regularity relative to the rank-one background $\varphi_n\otimes\varphi_n$ instead of the constant $\rho$. If the degrees in $V$ and in $V'$ each vary by a factor at most $1+\varepsilon$, then relative $\varepsilon^3$-regularity implies sparse $O(\varepsilon)$-regularity in the sense of Kohayakawa and R\"odl: $|d(S,T)-d(V,V')|=O(\varepsilon)\,d(V,V')$ whenever $|S|\ge\varepsilon|V|$ and $|T|\ge\varepsilon|V'|$.

\begin{theorem}[Degree partitions are regular]\label{thm:relreg}
Let $V_1,\dots,V_N$ partition $[n]$, $n_i=|V_i|$, $\mu_{ij}=\mu_{V_i,V_j}$, and $0<\varepsilon\le1$. If $\mu_{ij}\ge2\Delta_d/\varepsilon$, then
\[
\Pr\big[(V_i,V_j)\text{ is not $\varepsilon$-regular relative to }d\big]\le2\cdot4^{n_i+n_j}e^{-\varepsilon^2\mu_{ij}/2560}.
\]
\end{theorem}
\begin{proof}
Apply Theorem~\ref{thm:localstein} with $t=\varepsilon\mu_{ij}-\Delta_d\in[\varepsilon\mu_{ij}/2,\varepsilon\mu_{ij}]$, using $\mu_{S,T}\le\mu_{ij}$. Then $256(\mu_{S,T}+\Delta_d+t)\le640\mu_{ij}$, so both tails are at most $e^{-\varepsilon^2\mu_{ij}/2560}$. Take a union bound over $S\subseteq V_i$ and $T\subseteq V_j$.
\end{proof}

\begin{corollary}[A Szemer\'edi partition read off from the degrees]\label{cor:degpartition}
Let $\Lambda_n=\int\Psi(\varphi_n)\le\Lambda$, $M/n\to\infty$ and $K(d)=o(M)$ (for instance $\varphi_n\to\varphi$ in $L^1$ and $\rho\to0$). Fix $\varepsilon\in(0,\frac12]$ and $K=\exp(\Lambda/\varepsilon+1)$. Put the vertices with $\varphi_n<\varepsilon$ or $\varphi_n>K$ into an exceptional set $V_0$, and split the rest into the classes $\{\varepsilon(1+\varepsilon)^k\le\varphi_n<\varepsilon(1+\varepsilon)^{k+1}\}$. Then:
\begin{enumerate}[label=(\roman*)]
\item $V_0$ is incident to at most a fraction $2\varepsilon$ of the edge endpoints (by Lemma~\ref{lem:tail});
\item there are $N\le2(\Lambda/\varepsilon+1+\log\frac1\varepsilon)/\varepsilon+1$ classes;
\item with probability $1-o(1)$, every pair of classes with $\mu_{ij}\ge\varepsilon M/N^2$ is $\varepsilon$-regular relative to $d$, and the remaining pairs have total background mass $\sum\mu_{ij}\le\varepsilon M$ (hence, by Theorem~\ref{thm:stein}, at most $(\varepsilon+o(1))M$ edges).
\end{enumerate}
\end{corollary}
\begin{proof}
(i) and (ii) are immediate from the tail bound and from $\log(K/\varepsilon)/\log(1+\varepsilon)\le2\log(K/\varepsilon)/\varepsilon$. For (iii), $2\Delta_d/\varepsilon=o(M)$, so Theorem~\ref{thm:relreg} applies to all heavy pairs for large $n$, and the union bound over $N^2$ pairs gives $2N^2\,4^{n}e^{-\varepsilon^3M/(2560N^2)}\to0$, since $N$ is fixed and $M/n\to\infty$. The light pairs have $\sum\mu_{ij}<N^2\cdot\varepsilon M/N^2$.
\end{proof}

\begin{remark}[Meaning for the sparse regularity program]\label{rem:relregmeaning}
\begin{enumerate}[label=(\alph*)]
\item \emph{Size.} For uniformly random graphs with given degrees the regularity partition is explicit and has $O(\Lambda/\varepsilon^2+\varepsilon^{-1}\log\frac1\varepsilon)$ parts, against the generic $4^{O(1/\varepsilon^2)}$ for weak regularity (Theorem~\ref{thm:regularity}) and tower bounds for Szemer\'edi's lemma.
\item \emph{Energy equals entropy.} The only input from the profile is the modular $\Lambda$. By Lemma~\ref{lem:rankone}, $\Lambda(\varphi\otimes\varphi)\approx2\int\varphi\log\varphi$, which is also the entropy deficit per edge in $\log G(d)$. The same quantity bounds both the number of parts and the counting loss.
\item \emph{Relative regularity.} $G_d$ is quasirandom relative to its degrees, in the sense of Chung and Graham's discrepancy property for given degree sequences \cite{CG08}: there is no structure beyond the degree classes. This suggests that the right sparse regularity lemma is relative to a background, the degree-corrected rank-one graphon here and the $\beta$-graphon in the dense regime, in the way the Green--Tao transference is relative to a majorant.
\item \emph{Open.} Upgrade Theorem~\ref{thm:localstein} to a Bennett tail $\exp(-\mu\,\Psi(t/\mu))$, as Cook, Goldstein and Johnson obtain by size-biasing in the regular case \cite{CGJ18}. This would give the $L\log L$-shaped tail for hub blocks, in which $\Psi$ appears both as the regularity energy and as the concentration rate.
\end{enumerate}
\end{remark}

\subsection{One limit at all densities: the $\beta$-graphon, which is $\varphi\otimes\varphi$ when sparse}\label{subsec:betaphi}

The correct limit object at every density is the maximum-entropy $\beta$-model $w^\star$. In the sparse regime it coincides with $\varphi_n\otimes\varphi_n$ in the cut metric. The proof below is probabilistic: it applies Theorem~\ref{thm:stein} inside the $\beta$-model and uses neither the $\beta$-model equations nor any cap on $d_{\max}$.

\begin{theorem}[The $\beta$-model is $\varphi\otimes\varphi$ in the sparse sense]\label{thm:betaphi}
Let $d$ be strictly graphical with maximum-entropy probabilities $w^\star$, and suppose $n\ge6$ and $M\ge64n$. Then for every $L\ge2$,
\[
\big\|W^{w^\star}/\rho-\varphi_n\otimes\varphi_n\big\|_\square\le4\big(\mu_d(L)+\tau_d(L/2)\big)+20\sqrt{\frac nM}+\frac{5n}M .
\]
\end{theorem}
\begin{proof}
Let $G\sim G_{x^\star}$ and $d'=\deg(G)$, $M'=|E(G)|$, $D'_S=\sum_{i\in S}d'_i$. Conditionally on $d'$, $G$ is uniform with degree sequence $d'$. Since $E\,d'=d$ and $\Var d'_a=\sum_jw^\star_{aj}(1-w^\star_{aj})\le d_a$, we have $\sum_aE|d'_a-d_a|\le\sum_a\sqrt{d_a}\le\sqrt{2Mn}$. Fix $S,T$.

(1) $E\,e(S,T)=w^\star(S,T)$, because the edges are independent.

(2) Let $\mathcal E=\{M/2\le M'\le3M/2\}$. Bernstein's inequality gives $\Pr[\mathcal E^c]\le2e^{-3M/28}$, and on $\mathcal E^c$ all quantities are at most $n^2+2M$, so $\mathcal E^c$ contributes at most $1$ below.

(3) On $\mathcal E$, Theorem~\ref{thm:stein} for $d'$, integrated over $t$, gives
\[
E\big[|e(S,T)-D'_SD'_T/(2M')|\,\big|\,d'\big]\le4(2d'_{\max}+K(d'))+2\sqrt{8\pi M'}\le8n+4K(d')+13\sqrt M .
\]
By the splitting in \S\ref{subsec:stein}, $K(d')\le J'_L+\sum_{d'_a>L}d'_a$. Now $E J'_L\le J_L+\sum_aE|d'_a-d_a|$. A vertex with $d_a\le L/2<L<d'_a$ has $d'_a\le2|d'_a-d_a|$, so $E\sum_{d'_a>L}d'_a\le\sum_{d_a>L/2}d_a+2\sum_aE|d'_a-d_a|$. Hence
\[
E\,K(d')\le2M\big(\mu_d(L)+\tau_d(L/2)\big)+3\sqrt{2Mn}.
\]

(4) On $\mathcal E$, with $X_S=D'_S-D_S$ and $Z=M'-M$, we have $|D'_SD'_T/(2M')-D_SD_T/(2M)|\le2(|X_S|+|X_T|)+|X_SX_T|/M+4|Z|$. Since each edge changes $X_S$ by at most $2$, $\Var X_S\le4M$, so $E|X_S|\le2\sqrt M$, $E|X_SX_T|\le4M$ and $E|Z|\le\sqrt M$. This term is at most $12\sqrt M+4$.

Combining (1)--(4),
\[
\Big|w^\star(S,T)-\frac{D_SD_T}{2M}\Big|\le8M\big(\mu_d(L)+\tau_d(L/2)\big)+12\sqrt{2Mn}+25\sqrt M+8n+6 .
\]
Divide by $2M$ ($\int_{S\times T}W^{w^\star}/\rho=w^\star(S,T)/(2M)$), add the diagonal blocks of $\varphi_n\otimes\varphi_n$ (at most $n/(2M)$), and take the supremum over vertex sets; no union bound is needed.
\end{proof}

\begin{corollary}[Uniform graphs and the $\beta$-graphon]\label{cor:gdbeta}
Under the hypotheses of Theorem~\ref{thm:betaphi}, with probability at least $1-2e^{-n}$,
\[
\big\|W^{G_d}/\rho-W^{w^\star}/\rho\big\|_\square\le8\big(\mu_d(L)+\tau_d(L/2)\big)+25\sqrt{\frac nM}+\frac{10n}M .
\]
In particular, if $\varphi_n\to\varphi$ in $L^1$, $M/n\to\infty$ and $\rho\to0$, then $G_d$ and its $\beta$-model have the same normalized cut limit $\varphi\otimes\varphi$.
\end{corollary}
\begin{proof}
Theorems~\ref{thm:stein} and~\ref{thm:betaphi}, the triangle inequality, $\tau_d(L)\le\tau_d(L/2)$ and $d_{\max}\le n$.
\end{proof}

\begin{remark}[The uniform statement and what it still lacks]\label{rem:uniform}
Corollary~\ref{cor:gdbeta} is a single statement, $G_d\approx w^\star$, valid for every $d$, with an explicit remainder $\mu_d(L)+\tau_d(L/2)$: the share of degree mass that can form rich-club adjacencies, i.e.\ the \emph{rank-one defect} of the profile at scale $L$. Choosing $L=\sqrt M$, this remainder tends to $0$ for every uniformly integrable profile when $\rho\to0$, and then the $\beta$-graphon is $\varphi\otimes\varphi$. At fixed density it does not vanish, and correctly so, because the $\beta$-graphon then differs from $\varphi\otimes\varphi$ by order $\rho\iint\varphi^2\otimes\varphi^2$. At positive density, $G_d\approx w^\star$ with no remainder follows instead from a counting bound, $G(d)\ge\sqrt{D(d)}\,e^{-Cn\log n\log(2+M/n)}$ with $D(d)$ the number of loopless digraphs with in- and out-degrees $d$, which gives $-\log\pi_{x^\star}(d)=o(M)$ whenever $M\gg n\log n\log(M/n)$, and Theorem~\ref{thm:degconv} \cite{Meh-deg}. Numerically, for a profile $\varphi(q)\propto q^{-2/3}$ (tail exponent $2.5$) with $d_{\max}$ between $6\sqrt M$ and $11\sqrt M$, $\|W^{w^\star}/\rho-\varphi_n\otimes\varphi_n\|_1$ equals $0.201,0.191,0.184$ at $\rho=0.020,0.014,0.010$. This decreases slowly, as the tail mass $\int_{\varphi>\rho^{-1/2}}\varphi\asymp\rho^{1/4}$ predicts.
\end{remark}

\subsection{Entropy and $L\log L$}

\begin{lemma}\label{lem:rankone}
For measurable $\varphi\ge0$ with $\int\varphi=1$,
\[
2\int\varphi\log\varphi-1\;\le\;\Lambda(\varphi\otimes\varphi)\;\le\;2\int\varphi\log\varphi+2\log2.
\]
\end{lemma}
\begin{proof}
Lower bound: $\Psi(w)\ge w\log w-w$, and $\iint\varphi(x)\varphi(y)\log(\varphi(x)\varphi(y))=2\int\varphi\log\varphi$. Upper bound: $\Psi(w)\le w\log(1+w)$ and $\log(1+uv)\le\log(1+u)+\log(1+v)$ give $\Lambda(\varphi\otimes\varphi)\le2\int\varphi\log(1+\varphi)=2\int\varphi\log\varphi+2\int\varphi\log(1+1/\varphi)$, and by Jensen for the probability measure $\varphi\,dx$, $\int\varphi\log(1+1/\varphi)\le\log\int(\varphi+1)=\log2$.
\end{proof}

In the regime of Theorem~\ref{thm:degsparse}, the proof shows $\log G(d)=H_x(d)+O(n\log\bar d+d_{\max}^2)$, and expanding $H_x(d)$ gives
\[
\log G(d)=M\Big(\log\frac{n^2}{2M}+1-2\int\varphi_n\log\varphi_n\Big)+O\big(n\log\bar d+d_{\max}^2\big).
\]
So the same functional $\int\varphi\log\varphi$ measures the entropy deficit of a degree profile and, by Lemma~\ref{lem:rankone}, the $L\log L$ size of the limit $\varphi\otimes\varphi$ given by Theorem~\ref{thm:degsparse}.

\section{Sparse large deviations}\label{sec:ldp}

Extend $\Psi(u)=(1+u)\log(1+u)-u$ to $u\in[-1,\infty)$; it remains convex and nonnegative, with $\Psi''(u)=1/(1+u)$. Let
\[
h_p(x)=x\log\frac xp+(1-x)\log\frac{1-x}{1-p}\qquad(0\le x\le1)
\]
be the relative entropy of Bernoulli$(x)$ with respect to Bernoulli$(p)$; the rate of a graphon $W:[0,1]^2\to[0,1]$ for $G(n,p)$ is $I_p(W)=\frac12\int h_p(W)$ \cite{CV11}.

\begin{proposition}[$\Psi$ is the sparse rate]\label{prop:rate}
Let $0<p<e^{-2}$ and $x\in[0,1]$, and write $x=p(1+u)$. Then
\[
p\,\Psi(u)\;\le\;h_p(x)\;\le\;\big(1+\kappa(p)\big)\,p\,\Psi(u),\qquad\kappa(p)=2\sqrt p+\frac{2}{\log(1/p)-2}.
\]
Consequently, for every graphon $W$ with values in $[0,1]$ and $U=W/p-1$,
\[
\frac p2\int\Psi(U)\;\le\;I_p(W)\;\le\;\big(1+\kappa(p)\big)\frac p2\int\Psi(U).
\]
\end{proposition}
\begin{proof}
Let $D(x)=h_p(x)-p\Psi(x/p-1)$. Then $D(p)=D'(p)=0$ and $D''(x)=\frac1x+\frac1{1-x}-\frac1x=\frac1{1-x}\ge0$, so $D\ge0$ and $D(x)=\int_p^x\frac{x-t}{1-t}\,dt$. Hence $D(x)\le(x-p)^2$ for $x\le\frac12$ (as $\frac1{1-t}\le2$ between $p$ and $x$), and $D(x)\le x-p\le x$ for $x\ge p$ (as $\frac{x-t}{1-t}\le1$).

\emph{Case $x<p$.} Since $\Psi''\ge1$ on $[-1,0]$, $\Psi(u)\ge u^2/2$, so $p\Psi(u)\ge(x-p)^2/(2p)\ge D(x)/(2p)$, i.e.\ $D\le2p\cdot p\Psi(u)$.

\emph{Case $p\le x\le\sqrt p$.} Bernstein's inequality $\Psi(u)\ge u^2/(2(1+u/3))$ for $u\ge0$ \cite{BLM13} gives $p\Psi(u)\ge3(x-p)^2/(2(2p+x))$, so $D(x)\le(x-p)^2\le\frac23(2p+x)\,p\Psi(u)\le2\sqrt p\cdot p\Psi(u)$ (using $2p\le\sqrt p$).

\emph{Case $x\ge\sqrt p$.} $\Psi(u)\ge(1+u)\log(1+u)-(1+u)$, so $p\Psi(u)\ge x(\log\frac xp-1)\ge x(\frac12\log\frac1p-1)$, while $D(x)\le x$. Hence $D\le\frac{2}{\log(1/p)-2}\,p\Psi(u)$.

In all cases $0\le D\le\kappa(p)\,p\Psi(u)$. Integrate.
\end{proof}

\begin{example}[Clique and hub]\label{ex:cliquehub}
Let $a,b\in(0,1)$. If $W=1$ on $[0,a]^2$ and $W=p$ elsewhere (a planted clique on $an$ vertices), then $U=\frac1p-1$ on $[0,a]^2$ and $0$ elsewhere, and $\frac p2\Lambda(U)=\frac{a^2}2(\log\frac1p-1+p)$. If $W=1$ on $([0,b]\times[0,1])\cup([0,1]\times[0,b])$ and $p$ elsewhere (a hub of $bn$ vertices joined to everything), then $\frac p2\Lambda(U)=(b-\frac{b^2}2)(\log\frac1p-1+p)$. These are the two constructions behind sparse upper tails \cite{LZ17,BGLZ17}; the factor $\log\frac1p$ is $\Psi(K)/K$ at $K=1/p$, as in Lemma~\ref{lem:tail}.
\end{example}

\begin{remark}[No replica symmetry in the sparse regime]\label{rem:rsb}
For the upper tail $t(H,G)\ge(1+\delta)p^{e(H)}$ with $\Delta(H)\ge2$, the constant candidate $U\equiv(1+\delta)^{1/e(H)}-1$ has $\frac p2\Lambda(U)=\Theta(p)$, i.e.\ cost of order $n^2p$, whereas planted structures achieve the same increase at cost of order $n^2p^{\Delta}\log\frac1p=o(n^2p)$ \cite{BGLZ17}. So in the sparse regime the constant graphon is never optimal, and the question settled by generalized H\"older in the dense regime \cite{LZ15} (when is the constant optimal?) has no sparse analogue. In the sparse range $n^{-1/\Delta}\ll p\ll1$ the variational problem was solved by Lubetzky--Zhao for cliques \cite{LZ17} and by Bhattacharya--Ganguly--Lubetzky--Zhao in general \cite{BGLZ17}; for regular $H$ the optimum is the cheaper of a clique and a hub, described through an independence polynomial, and for irregular $H$ a hub construction is optimal.
\end{remark}

\begin{remark}[Mean-field theorems]\label{rem:meanfield}
Nonlinear large deviation theory \cite{CD16,Eld18,CoDe20,Aug20,HMS22} shows, in increasingly wide ranges of $p$, that $-\log\Pr[t(H,G(n,p))\ge(1+\delta)p^{e(H)}]$ is asymptotic to the value of the $n$-vertex discrete variational problem $\inf\{\sum_{i<j}h_p(q_{ij}):\E_q\,t(H,G_q)\ge(1+\delta)p^{e(H)}\}$, and its asymptotics in the range above are given by \cite{BGLZ17}; see \cite{BB23} for the localized regime. By Proposition~\ref{prop:rate}, each term $h_p(q_{ij})$ may be replaced by $p\Psi(q_{ij}/p-1)$ at the cost of a factor $1+O(1/\log\frac1p)$, so the sparse upper-tail rate is, to leading order, a minimization of the modular $\Lambda$ under the truncation $U\le\frac1p-1$.
\end{remark}

\begin{remark}[Why compactness alone does not give the LDP]
The Chatterjee--Varadhan route \cite{CV11} derives the dense LDP from compactness of graphon space together with continuity of $t(H,\cdot)$ in the cut metric. Compactness survives in $L\log L$ (Theorem~\ref{thm:compact}), but continuity does not: $t(H,\cdot)$ can be infinite on $L\log L$ balls (Section~\ref{sec:open}), so no bound of the form $t(H,W)\le f(\Lambda(W))$ exists. The truncation $U\le\frac1p-1$ is what makes the sparse problem well posed.
\end{remark}

\begin{proposition}[Transfer to degree-constrained graphs]\label{prop:transfer}
Let $G_\beta$ be a $\beta$-model, i.e.\ a random graph with independent edges of probabilities $w_{ij}=x_ix_j/(1+x_ix_j)$, let $d$ be a degree sequence with $\pi(d)=\Pr[\deg(G_\beta)=d]>0$, and let $G_d$ be uniformly random with degree sequence $d$. For every event $\mathcal A$,
\[
\Pr[G_d\in\mathcal A]\;\le\;\Pr[G_\beta\in\mathcal A]/\pi(d).
\]
\end{proposition}
\begin{proof}
Every graph with degree sequence $d$ has probability $\prod_{i<j}(1-w_{ij})\prod_ix_i^{d_i}$ under $G_\beta$, so $G_\beta$ conditioned on $\deg(G_\beta)=d$ is uniform \cite{CDS11}. Hence $\Pr[G_d\in\mathcal A]=\Pr[G_\beta\in\mathcal A,\deg=d]/\pi(d)$.
\end{proof}

Thus upper bounds for upper tails transfer from the product measure $G_\beta$, to which nonlinear large deviation theory applies, to uniformly random graphs with given degrees, at the cost $-\log\pi(d)$. When $x$ is the maximum-entropy choice, so that $\E\deg(G_\beta)=d$, this cost is $O(n\log(2+M/n))$ for sparse sequences with $d_{\max}=o(M^{1/4})$, by the McKay--Wormald asymptotics \cite{McK85,MW91}, and for dense sequences whose maximum-entropy edge probabilities are bounded away from $0$ and $1$, by Barvinok--Hartigan \cite{BH13}. This is negligible against upper-tail rates of order $n^2p^\Delta\log\frac1p$ whenever $np^\Delta\log\frac1p\gg\log(2+M/n)$.

\section{Open problems}\label{sec:open}

\begin{enumerate}[leftmargin=2em,label=\arabic*.]
\item \emph{Counting.} For a graph $F$ of maximum degree $\Delta\ge2$, $t(F,W)$ can be infinite for $W\in L\log L$: for the path with two edges, $t(P_2,\varphi\otimes\varphi)=\int\varphi^2\cdot(\int\varphi)^2$, which is infinite in Example~\ref{ex:beyond}. In \cite{BCCZ19} counting lemmas require $p>\Delta$ and fail for $p\le\Delta$ \cite[Prop.~2.22]{BCCZ19}. Determine which normalized or truncated homomorphism densities are continuous on $L\log L$ balls.
\item \emph{Vertex partitions.} Prove the equivalence of the graphon-partition and vertex-partition definitions of upper $\Psi$-regularity for simple graphs with average degree tending to infinity.
\item \emph{The bound on $\eta$.} Is the factor $e^{-O(C/\varepsilon)}$ in the admissible $\eta$ of Theorem~\ref{thm:regularity} necessary?
\item \emph{Signed graphons.} Extend the theory to signed graphons and weighted graphs with no dominant nodes, as in \cite{BCCZ19}.
\item \emph{Degree-constrained graphs: rates and local statistics.} Convergence of $G_d$ to its $\beta$-model holds at every density with $M\gg n\log n\log(M/n)$ \cite{Meh-deg}, and in the sparse regime for every $L^1$ profile (Corollary~\ref{cor:stein}). Find optimal rates in the cut metric, extend Theorem~\ref{thm:localstein} to a Bennett tail $\exp(-\mu\Psi(t/\mu))$, and identify which local statistics converge; see \cite{Meh-open}.
\item \emph{Large deviations via Young's inequality.} Prove the lower bound in the sparse upper-tail variational problem of \cite{BGLZ17}, in the modular form of Remark~\ref{rem:meanfield}, using only the Young pair $(\Psi,\Psi^*)$ and the truncation $U\le\frac1p-1$.
\item \emph{Upper tails for degree-constrained graphs.} Using Proposition~\ref{prop:transfer}, derive upper-tail bounds for subgraph counts in uniformly random graphs with a given degree sequence, and identify matching lower bounds.
\end{enumerate}

\end{document}